\documentclass[12pt]{article}
\usepackage{epsfig}
\usepackage{color}
\usepackage{amsfonts,amsthm,amsmath,amssymb,amscd,graphicx}

\begin{document}
\newtheorem{cor}{Corollary}
\newtheorem{lm}{Lemma}
\newtheorem{theorem}{Theorem}
\newtheorem{df}{Definition}
\newtheorem{prop}{Proposition}
\newtheorem{remark}{Remark}

\title{On a scenario of onset of strongly dissipative mixed dynamics}
\author{Kazakov A.O.$^{1}$\\
$^1$National Research University Higher School of Economics, \\
25/12 Bolshaya Pecherskaya Ulitsa, 603155 Nizhny Novgorod, Russia.
}
\maketitle

\begin{abstract}
In this paper we present the scenario of the occurrence of strongly dissipative mixed dynamics in two-dimensional reversible diffeomorphisms, using as an example the system describing a  motion of two point vortices under the influence of wave perturbation and shear flow. For mixed dynamics of this type the chaotic attractor intersects with the chaotic repeller, but their intersection forms a ``thin'' set. The main stage of this scenario is the appearance of homoclinic structures for a symmetric saddle orbit which arise after crisis of a homoclinic attractor and repeller.
\end{abstract}

\section{Introduction}

{\bf On three types of chaos}

Until recently, it was believed that chaos in dynamical systems can be only of two forms: conservative and dissipative chaos. Conservative (Hamiltonian) chaos is typical for systems preserving the phase volume, while, the dissipative chaos is often observed in systems in which the phase volume is compressed. It follows from the well-known theorem of Conley \cite{Conley78} that in any dynamical system defined on a compact manifold there exists an attractor and a repeller. There are many different definitions of an attractor. Here, as usual, by an attractor we will mean a stable closed invariant set of a system, and be a repeller -- an attractor under the time reversal. However, following Conley \cite{Conley78}, Ruelle \cite{Ruelle1981}, and Hurley \cite{Hurley1982} under stability we will mean the so-called stability under permanent acting perturbations, which is also called {\it total stability} or {\it Lyapunov stability by $\varepsilon$-orbits}, see more details in \cite{GonTur2017}.

We note, that such defined attractors (repellers) are in a good agreement with numerics and, what is important, allow to distinguish three types of chaos: conservative, dissipative and {\it mixed dynamics} \cite{GonTur2017}. So, if we denote the attractor in the system by $\mathcal A$ and the repeller by $\mathcal R$ then the condition $\mathcal A = \mathcal R$ holds for conservative chaos (moreover in this case, the attractor / repeller is the whole phase space), while the condition $\mathcal A \cap \mathcal R = \emptyset$ takes always place for dissipative chaos. The mixed dynamics is characterized by the fact, that the attractor and the repeller intersect, but do not coincide \cite{GonTur2017}, i.e. $\mathcal A \cap \mathcal R \neq \emptyset$, and $\mathcal A \neq \mathcal R$. It is worth noticing that, from a logical point of view, this is the third and final possibility and,  thus, a fourth type chaos does not exists.



For the first time the term ``mixed dynamics'' for this phenomenon was introduced in \cite{GonShilSten2002} (and later was used, e.g. in the papers, \cite{GonShilSten2006, GonDelsh2012}). However, the phenomenon associated with the possibility of intersecting attractors and repellers was discovered earlier in the paper \cite{GonShilTur97} where it was proved that near any two-dimensional diffeomorphism with a non-transverse heteroclinic cycle containing 2 saddle points with Jacobians less and greater than one\footnote{Here, under non-transverse heteroclinic cycle we mean a contour containing 2 saddle periodic points for which one pair of stable and unstable invariant manifolds intersects transversally and the other pair have a heteroclinic tangency.} there exist open regions (Newhouse domains) in which there are dense diffeomorphisms with infinitely many stable, completely unstable, and saddle periodic orbits whose closures have non-empty intersections, which, in fact, means that attractors are not separable from repellers. Later, in the paper \cite{LambStenkin2004}, an analogous theorem was proved in the case when a diffeomorphism $f$ with such a non-transverse heteroclinic cycle is {\it reversible}, i.e., when $f$ and $f^{-1}$ are conjugated by means of some involution $h$ for which $h \circ h = id$. However, in contrast to the general case, {\it symmetric periodic orbits} of conservative type (with the Jacobian $J = 1$), such as elliptic and area-preserving saddles, also arise here in addition to periodic sinks and sources. Note, that such symmetric periodic orbits intersect the set $Fix(h)$ of fixed points of the involution $h$ (i.e. the set of such points $x$ that $h(x) = x$). Such orbits appear due to the {\it reversible Newhouse phenomenon} \cite{LambStenkin2004}.

{\bf Mixed dynamics in reversible systems}

The phenomenon, when the dissipative dynamics in the system coexist with the conservative dynamics, was previously discovered in the papers \cite{PolitiOppoBadii86, QuispelRoberts92, LambRoberts98}, where it was found that the phase space for reversible two-dimensional diffeomorphisms can be divided into invariant domains with conservative dynamics\footnote{Probably, inside such invariant domains there exist inseparable sets of stable and completely unstable periodic orbits. But the dynamics in such areas looks conservative due to the fact that the absolute values of Jacobians for such orbits are very close to 1.} and invariant domains containing pairs of attractor--repeller. The possibility of intersecting an attractor and a repeller in physical systems was (as we know) first found numericaly in \cite{PikTop2002}, where for the model of four coupled rotators was observed that the attractor and the repeller ``may overlap''. Later, in \cite{GonGonKazTur2017}  such ``overlaping'' was explained by the emergence of the mixed dynamics in this model. Moreover, in \cite{GonGonKazTur2017} two different scenarios of the occurrence of mixed dynamics were also described. According to the first scenario, mixed dynamics appear in a ``soft'' manner, due to a sequence of local and global bifurcations of symmetry breaking. According to the second scenario, mixed dynamics can arise by explosion due to the emergence of heteroclinic structures, which appear after crisis of a simple attractor and a simple repeller.

Among the studies of mixed dynamics in systems from applications, here we would like to note the papers \cite{GonGonKaz2013, Kaz2013, Kazakov2015} in which mixed dynamics were found in nonholonomic models of the Celtic stone, rubber Chaplygin top, and Suslov top respectively.

\begin{figure}[!ht]
\begin{minipage}[h]{0.49\linewidth}
\center{\includegraphics[width=1\linewidth]{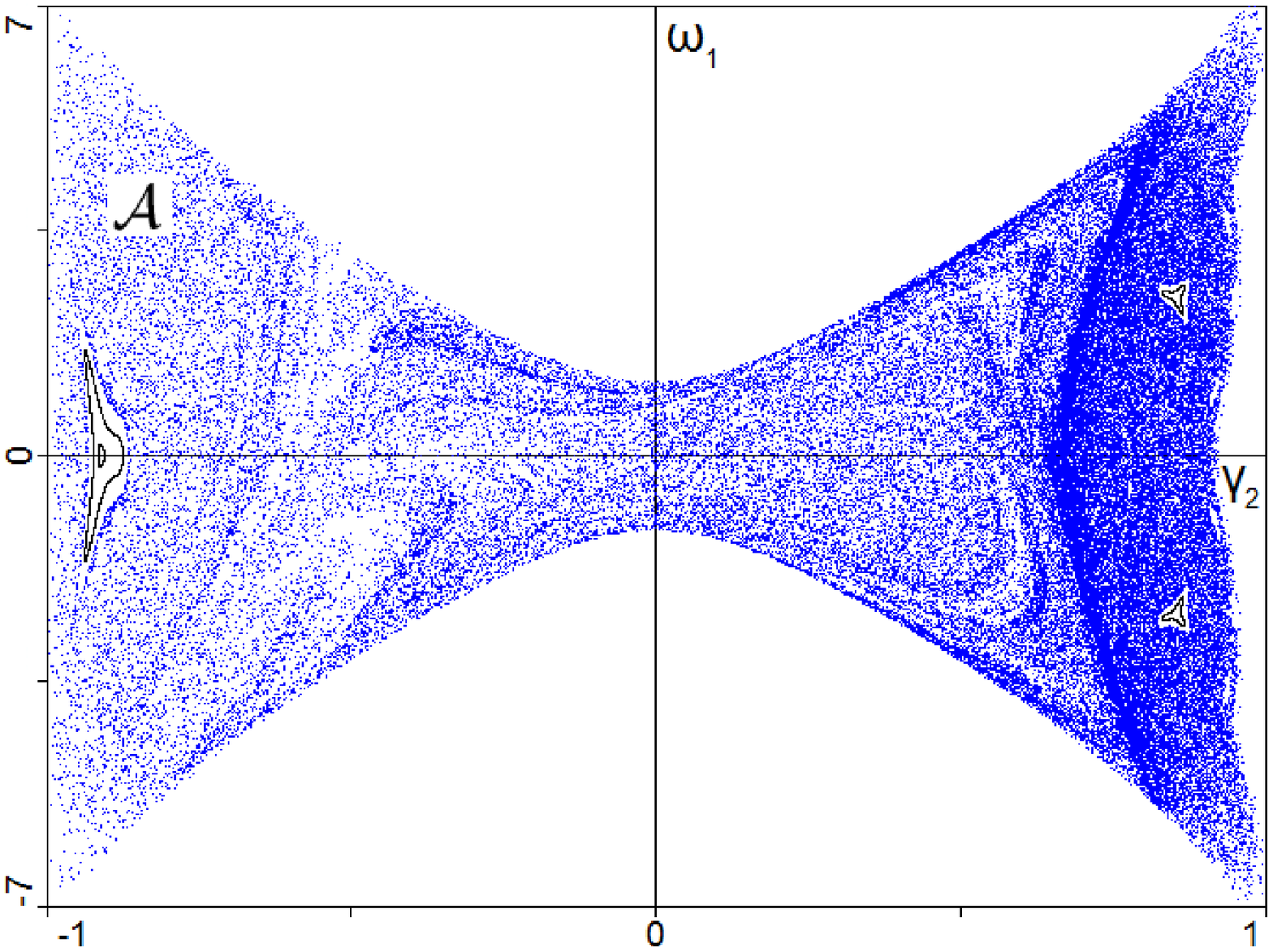} \\ a) {\footnotesize Attractor $\mathcal A$}}
\end{minipage}
\hfill
\begin{minipage}[h]{0.49\linewidth}
\center{\includegraphics[width=1\linewidth]{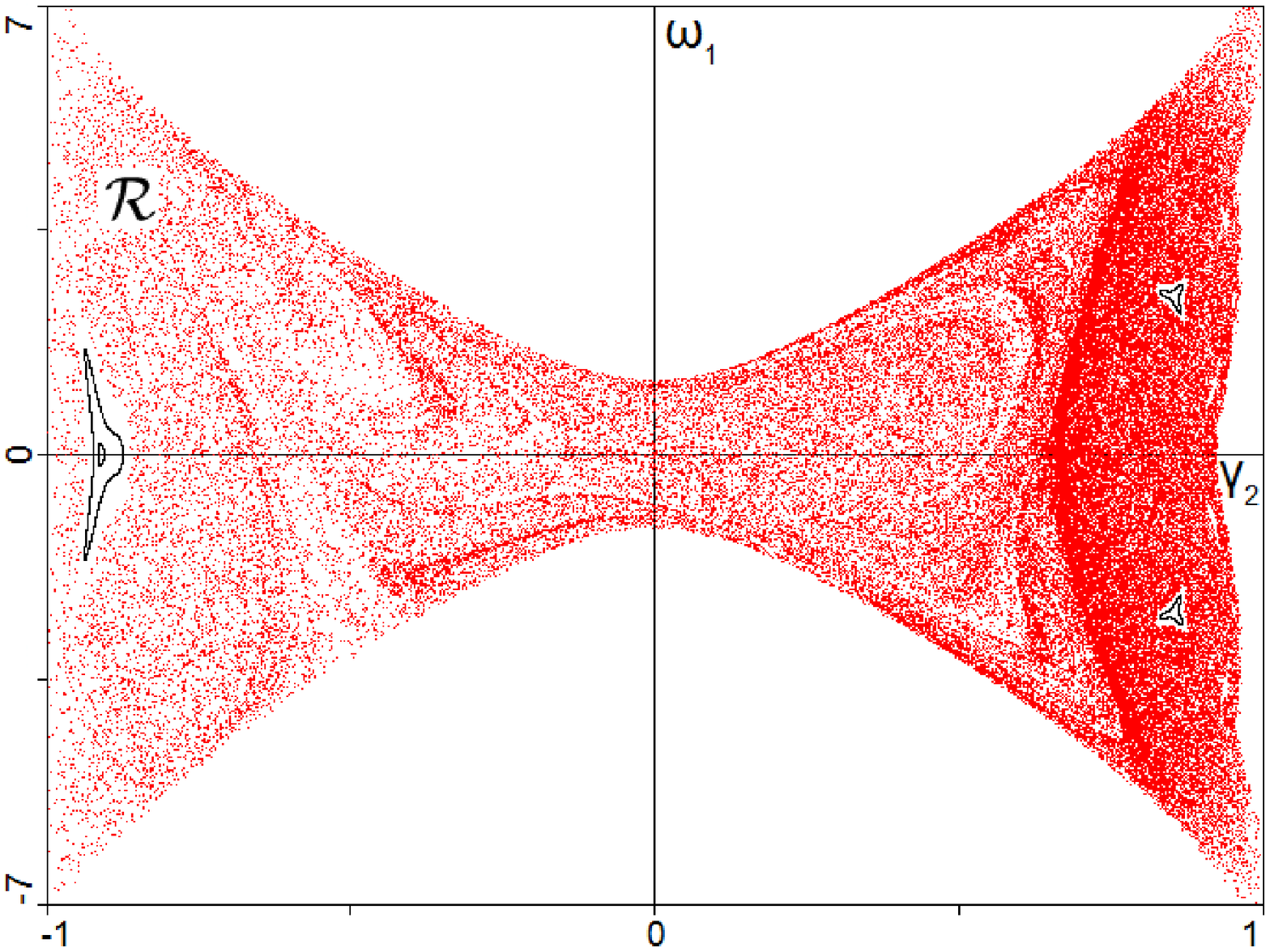} \\ b) {\footnotesize Repeller $\mathcal R$}}
\end{minipage}
\caption{{\footnotesize The phase portraits of the attractor $\mathcal A$ and the repeller $\mathcal R$ in the nonholonomic model of Suslov top \cite{BizBorKaz2015}. The line $\omega_1 = 0$ forms the set of fixed points of the involution $h : \omega_1 \rightarrow -\omega_1$ in this model. It is clearly seen
that $\mathcal A \cap \mathcal R \neq \emptyset, \mathcal A \neq \mathcal R$.}}
\label{fig:MD_Suslov}
\end{figure}

We note that in almost all cases listed above, mixed dynamics appears in reversible systems, which are obtained as a result of small perturbation of conservative systems. As a result of such perturbations the conservativity is destroyed, in general. Thus, the phase portraits of attractors and repellers will be slightly different but, due to reversibility, they remains to be symmetric with respect to $Fix(h)$ (see e.g. Fig. \ref{fig:MD_Suslov}, where the phase portraits of the attractor and the repeller in the nonholonomic model of Suslov top \cite{Kazakov2015} are presented).



{\bf Strongly dissipative mixed dynamics}

\begin{figure}[!ht]
\begin{minipage}[h]{0.49\linewidth}
\center{\includegraphics[width=1\linewidth]{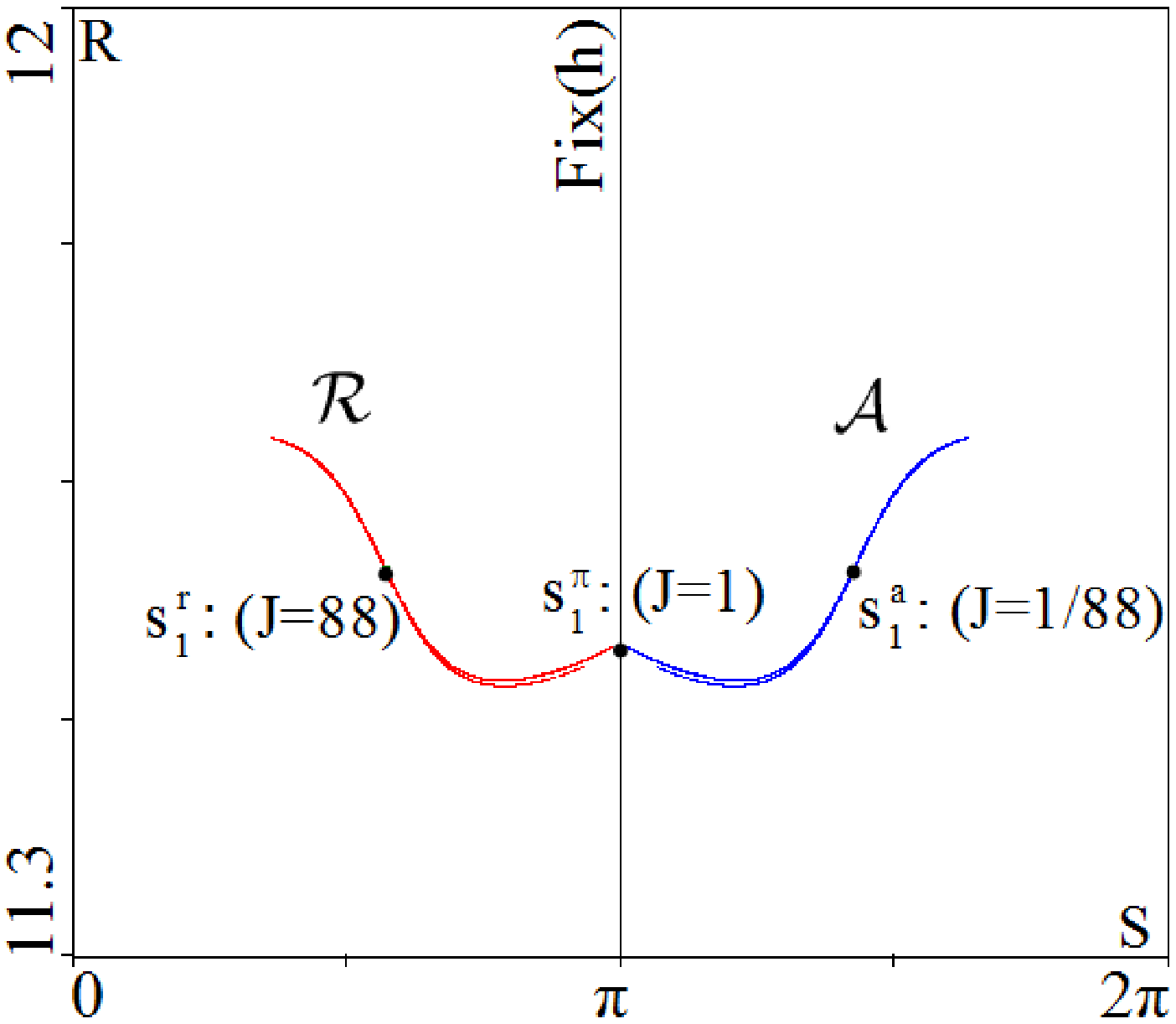} \\ a) {\footnotesize $\varepsilon = 0.1463$}}
\end{minipage}
\hfill
\begin{minipage}[h]{0.49\linewidth}
\center{\includegraphics[width=1\linewidth]{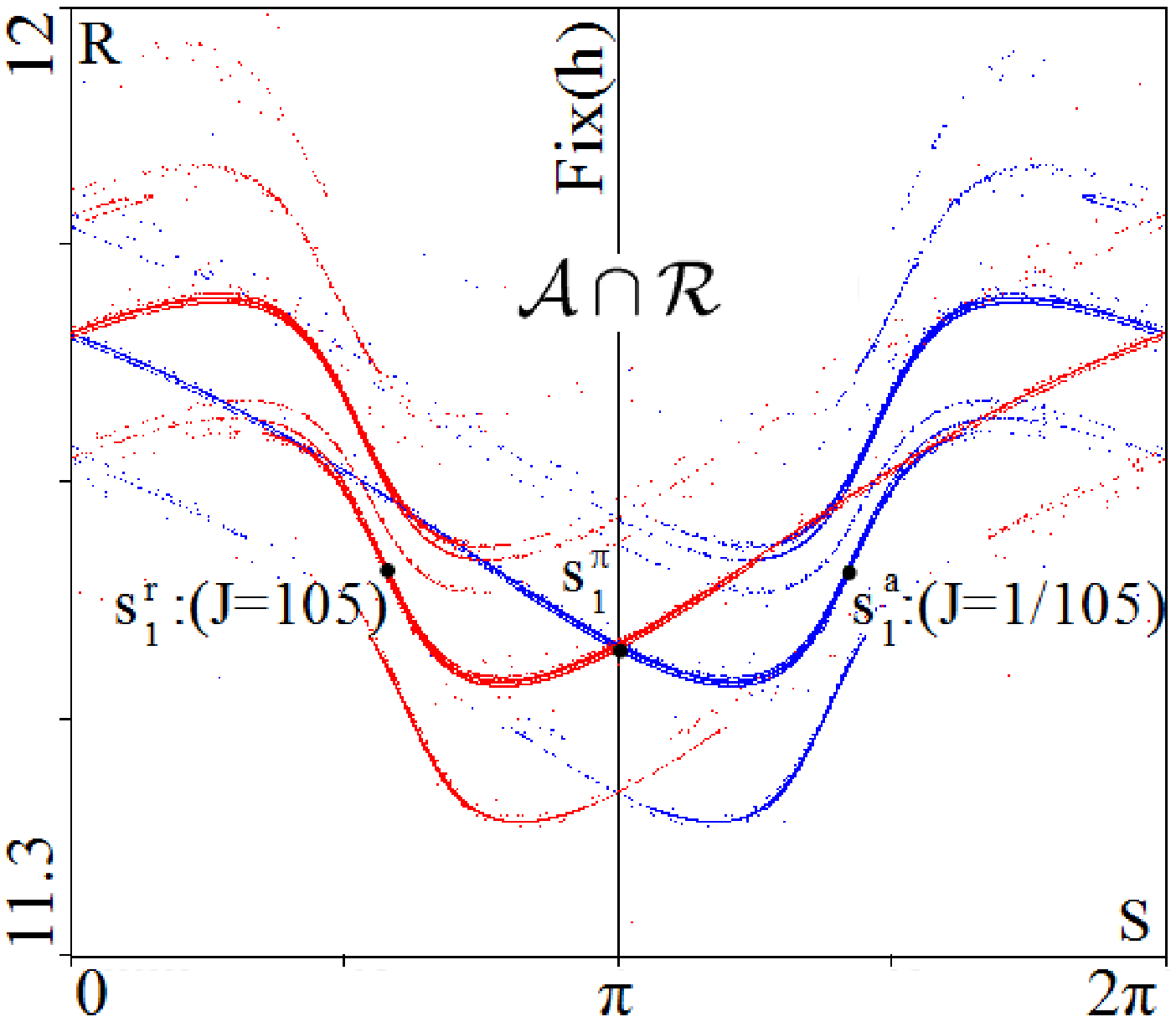} \\ b) {\footnotesize $\varepsilon = 0.14815$}}
\end{minipage}
\caption{{\footnotesize The phase portraits of the attractor $\mathcal A$ and the repeller $\mathcal R$ in the model of two point vortices perturbed by a wave and a shear flow. Fix(h) -- is the line of fixed points of the involution $h: \{S \rightarrow 2\pi - S, \quad R \rightarrow R\}$, $s_1$ -- is a symmetrical saddle point (with Jacobian $J = 1$). (a) $\mathcal A \cap \mathcal R = \emptyset$. Saddle point $s_a$ with Jacobian $J = 1/88$ belongs to the attractor, while the saddle point $s_r$ with $J = 88$ belongs to the repeller. (b) $\mathcal A \cap \mathcal R \neq \emptyset, \mathcal A \neq \mathcal R$. Mixed dynamics contain dissipative and also conservative periodic orbits.}}
\label{fig:MD_Vortex}
\end{figure}

In this paper we present the scenario of the appearance of strongly dissipative mixed dynamics in reversible two-dimensional diffeomorphisms, using as an example the system describing the motions of two point vortices perturbed by a wave and a shear flow \cite{Vet2017}. For mixed dynamics of this type the attractor and the repeller are very different from each other, and their intersection forms a ``thin'' set (see. Fig. \ref{fig:MD_Vortex}b). The main stage of this scenario is the appearance of homoclinic structures (for a symmetrical saddle orbit $s_1$) which arise after the crisis of a homoclinic attrator and repeller. Recall, that under the homoclinic attractor here we mean a chaotic attractor which contains a fixed point of a saddle type \cite{GonGon2016}. In the case under consideration, the attractor $\mathcal A$ contains a saddle point with very small Jacobian $J$ and the repeller $\mathcal R$ -- with very large Jacobian $J^{-1}$ (for example, in the case presented in Fig. \ref{fig:MD_Vortex}a, $J = 1/88$ for the attractor, while $J^{-1} = 88$ for the repeller). After the appearance of mixed dynamics, both the attractor and the repeller also contain a symmetrical saddle point $s_1$ with $J = 1$ (see. Fig. \ref{fig:MD_Vortex}b).

We note that, in accordance with Theorem 1 from \cite {GonTur2017}, the set belonging to the the intersection of an attractor and a repeller (also called the {\it reversible core}) is the limit of an infinite sequence of attractors, and also of an infinite sequence of repellers. Therefore, when the saddle point $s_1$ belongs to the reversible core, it can be argued that strongly dissipative mixed dynamics also contain countable sets of stable, as well as completely unstable periodic orbits.

However, in contrast to previously investigated cases (\cite{GonGonKaz2013, Kaz2013, Kazakov2015, GonGonKazTur2017}) where stable, as well as completely unstable periodic orbits, belonging to the regions with mixed dynamics, have Jacobians close to 1, strongly dissipative mixed dynamics contains periodic saddles, sinks and sources with various Jacobians (large, small, close, and equal to 1). In our opinion, the existence of such dissipative orbits is explained by the fact that the system strongly contracts volumes in some regions of phase space and strongly expands volumes in other (symmetrical) regions. We also suppose, that such a property of the system leads to a large differences between numerically obtained attractor and repeller after their intersection (when mixed dynamics appears)\footnote{We note, that in this case, theoretically, the attractor should be slightly different from the repeller due to appearance of ``homoclinic tangle'' \cite{GonTur2017}. But due to finiteness of the numerical calculation, the difference between attractor and repeller remains noticeable.}.

The work is organized as follows: the scenario of the appearance of strongly dissipative mixed dynamics in two-dimensional reversible diffeomorphisms is described in Sec. \ref{sec:Scenario}, the example of the implementation of such a scenario in the equations describing the motions of two point vortices perturbed by a wave and a shear flow is given in Sec. \ref{sec:VortexScenario}, and the equations governing the motions in this model are presented in Sec. \ref{sec:Eq}.

\section{Scenario of the appearance of strongly dissipative mixed dynamics in reversible systems} \label{sec:Scenario}

\begin{figure}[!ht]
\center{\includegraphics[width=0.99\linewidth]{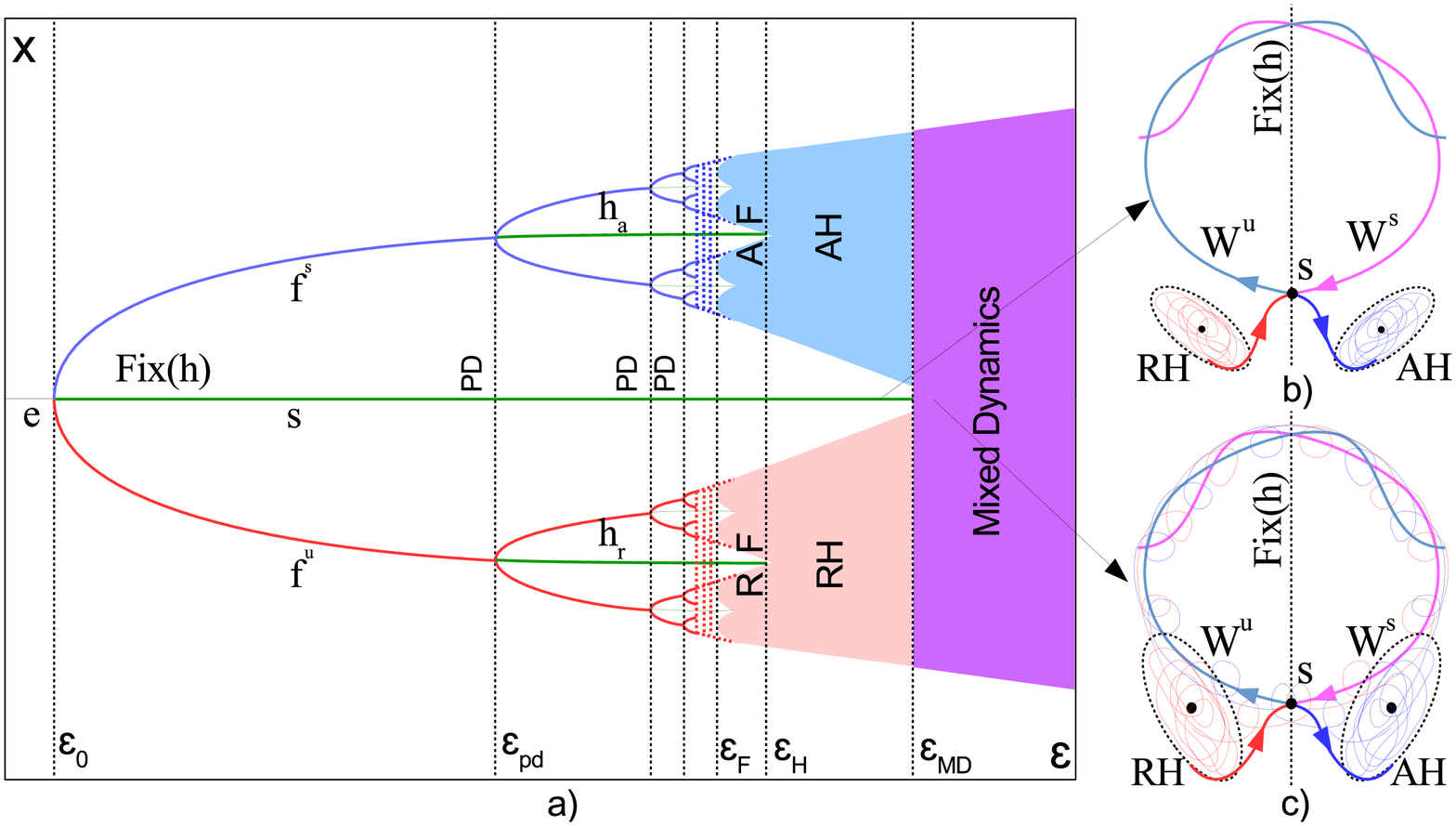} }
\caption{{\footnotesize The scenario of the appearance of strongly dissipative mixed dynamics in reversible two-dimensional diffeomorphisms.}}
\label{fig:MD_scenario_scheme}
\end{figure}

Let us consider a one-parameter family of two-dimensional reversible diffeomorphisms $\bar X = F(X,\varepsilon)$ defined on a compact manifold and depending on the parameter $\varepsilon$.
For simplicity we suppose that this family has a symmetric elliptic fixed point $e$ for $\varepsilon < \varepsilon_0$ and this point (at $\varepsilon = \varepsilon_0$) undergoes reversible pitch-fork bifurcation \cite{LermanTuraev2012}  after which this point becomes a symmetric saddle point $s$ and a pair of stable $f^s$ and completely unstable $f^u$ points (one point is symmetric to another) appears in a neighborhood of $s$ (see Fig. \ref{fig:MD_scenario_scheme}, near $\varepsilon = \varepsilon_0$). We also suppose that, with further increase of $\varepsilon$, a {\it Feigenbaum-like attractor} $AF$ \cite{Feigenbaum83} is born as a result of the cascade of period-doubling bifurcations of the point $f^s$. In the same way, a Feigenbaum-like repeller $RF$ is born from $f^u$. We note, that after the first period doubling bifurcation the points $f^s$ and $f^u$ become saddle: we denote these fixed points as $h_a$ and $h_r$.

Recall that immediately after the onset of chaotic dynamics through the cascade of period-doubling bifurcations, the Feigenbaum-like attractor consists of a set of components. With the change of the parameter the components of this attractor pairwise merge (as a result of the occurrence of heteroclinic intersections between manifolds of the saddle orbits belonging to the components of the attractor and manifolds of the saddle orbits lying between these components). Finally, two last components separated by $h_a$ are merged and the {\it homoclinic Henon-like attractor} $AH$ appears \cite{Henon76} (see Fig. \ref{fig:MD_scenario_scheme}a at $\varepsilon = \varepsilon_H$ and Fig. \ref{fig:MD_scenario_scheme}b). Symmetrically the homoclinic Henon-like repeller $RH$ containing the fixed point $h_r$ occurs.

With a further increase in $\varepsilon$, the Henon-like attractor $AH$ becomes larger and approaches the boundary of its basin of attraction which is formed by the stable manifold $W^s$ of the saddle point $s$ (accordingly, the basin for the Henon-like repeller $RH$ is bounded by the unstable manifold $W^u$ of the same point $s$). Also we note that starting with a certain value of the parameter the pairs of manifolds $W^s$ and $W^u$ intersects transversally, i.e. formed the homoclinic structure. (see Fig. \ref{fig:MD_scenario_scheme}b).

When $\varepsilon = \varepsilon_{MD}$, the crisis of the attractor $AH$ and the repeller $RH$ appears ($AH$ collides with the boundary of its basin of attraction $W^s$, while $RH$ collides with $W^u$), after which both these sets capture the homoclinic structure of the saddle point $s$ (see Fig. \ref{fig:MD_scenario_scheme}b)\footnote{Similarly, the crisis of the homoclinic attractor in Henon map appears \cite{Henon76}. However, in the case of Henon map, after the intersection of the attractor with the stable manifold bounding its basin of attraction, almost all trajectories go to infinity.}. In the case when trajectories in the system do not tend (in forward time) to other attractors\footnote{When the system is multistable, trajectories from the neighborhood of just collapsed attractor can tend to other attractors that coexist with $AH$. In this case mixed dynamics exist but do not manifest itself.}, since both the attractor and the repeller contain this homoclinic structure, it can be stated that the attractor intersects with the repeller, i.e. mixed dynamics appears.

\section{The model under consideration}\label{sec:Eq}

We consider a system describing the motion of two identical point vortices under the influence of a wave perturbation and external shear flow with a uniform distribution of vorticity \cite{Vet2017} as an example of a model in which the scenario described above is realized and strongly dissipative mixed dynamics appears.
\begin{equation}
\begin{cases}
\dot R = \frac{1}{2} A R \sin{2\varphi} - \varepsilon \sin{\varphi}\sin{S}\sin(R \sin{\varphi}) \\
\dot S = -1 + \varepsilon \cos{S} \cos(R \sin{\varphi}) \\
\dot \varphi = \frac{\kappa}{R^2} + A \cos^2{\varphi} - \frac{\varepsilon}{R} \cos{\varphi} \sin{S} \sin(R \sin{\varphi}).
\end{cases}
\label{eq:mainEq}
\end{equation}
Here $R \in (0, \infty), S \in [0, 2\pi), \varphi \in [0, 2\pi)$ are the phase variables and $A, \varepsilon, \kappa$ are the parameters of the system. Note that this model is a generalization of the system describing the motions of two point vortices which interact with a potential wave \cite{GoncharOstapchukTur91, KazVet2016}\footnote{The author of the system is E.V. Vetchanin \cite{Vet2017}. In contrast to the system investigated in \cite{GoncharOstapchukTur91, KazVet2016}, here the vortices also perturbed by an external flow with constant vorticitiy $A$ in addition to the interaction with a potential wave.}. In this system the parameter $\varepsilon$ specifies the amplitude of a wave, $A$ is the vorticity of the external flow, and $\kappa$ is the sum of intensities of the vortices.

Note that the equations \eqref{eq:mainEq} are invariant with respect to the substitution
\begin{equation}
\begin{gathered}
H_1: \{R \rightarrow R,\quad S \rightarrow -S, \quad \varphi \rightarrow -\varphi, \quad t \rightarrow -t\},\\
H_2: \{R \rightarrow R,\quad S \rightarrow 2\pi - S, \quad \varphi \rightarrow 2\pi - \varphi, \quad t\rightarrow -t\}.
\label{eq:involFlow}
\end{gathered}
\end{equation}
Thus, the system under consideration is reversible and $H = H_1 \cup H_2$ is the involution of this system. Also we note, that the variable $\varphi$ can be chosen as a secant in the system. In the case $\varphi = 0$ (or $\varphi = \pi$) the involution $H$ defines the involution $h = h_1 \cup h_2$ for the Poincar\'e map of the system, where:
\begin{equation}
\begin{gathered}
h_1: \{R \rightarrow R,\quad S \rightarrow -S\}\\
h_2: \{R \rightarrow R,\quad S \rightarrow 2\pi - S\},
\label{eq:involMap}
\end{gathered}
\end{equation}
and the set $Fix(h)$ of fixed points of this involution consists of two lines:
$$
Fix(h) = \{S = 0\} \cup \{S = \pi\}.
$$

Hereafter we choose the variable $\varphi = 0 $ as a secant for the system and perform one-parameter analysis by varying $\varepsilon$, assuming that other parameters are fixed as follows:
$$
A = 0.1, \kappa = 4.65.
$$

\section{A description of a scenario of occurrence of strongly dissipative mixed dynamics in the vortex model under consideration}\label{sec:VortexScenario}

\begin{figure}[!ht]
\begin{minipage}[h]{0.49\linewidth}
\center{\includegraphics[width=1\linewidth]{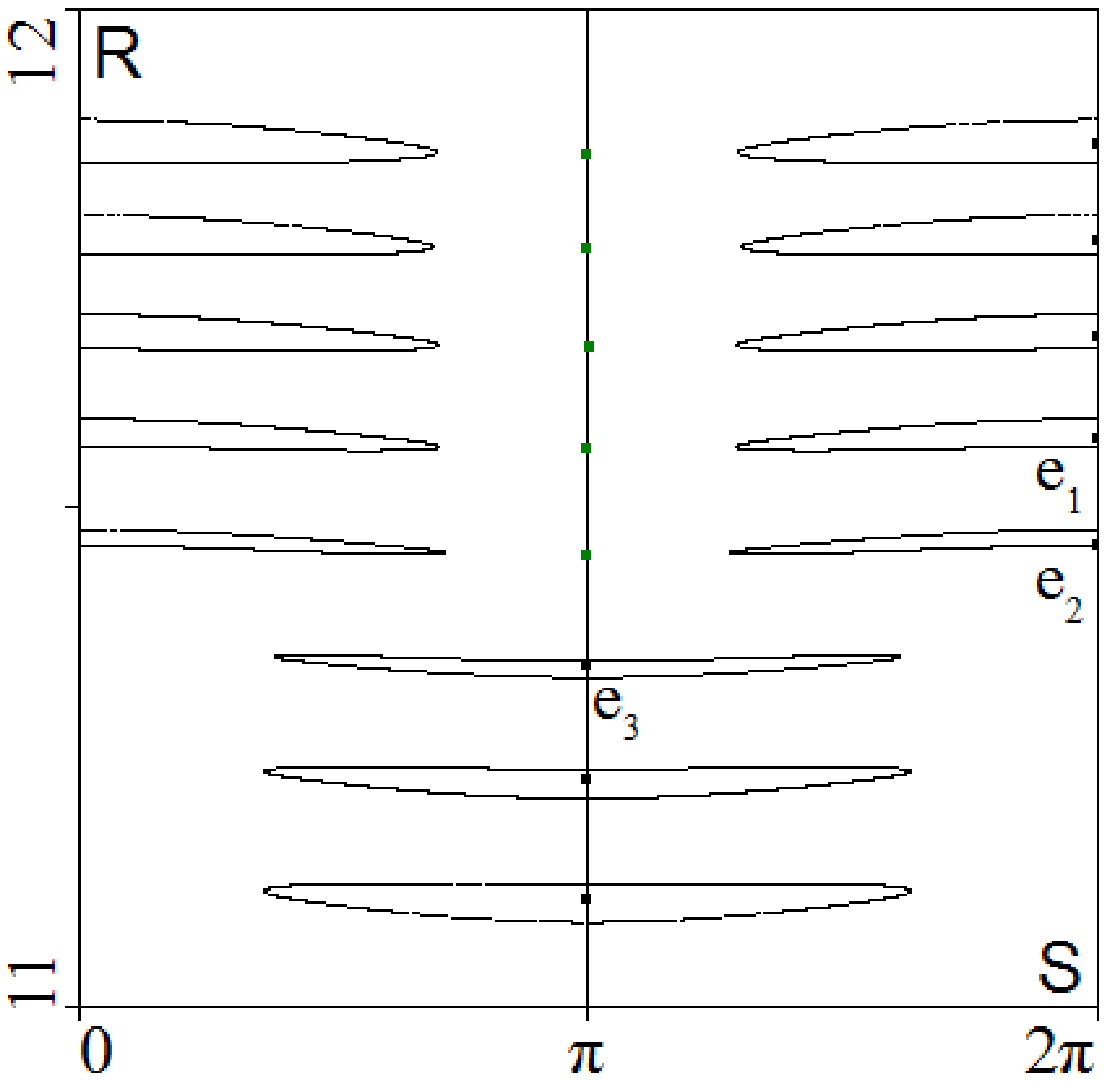} \\ a) {\footnotesize $\varepsilon = 0.01$}}
\end{minipage}
\hfill
\begin{minipage}[h]{0.49\linewidth}
\center{\includegraphics[width=1\linewidth]{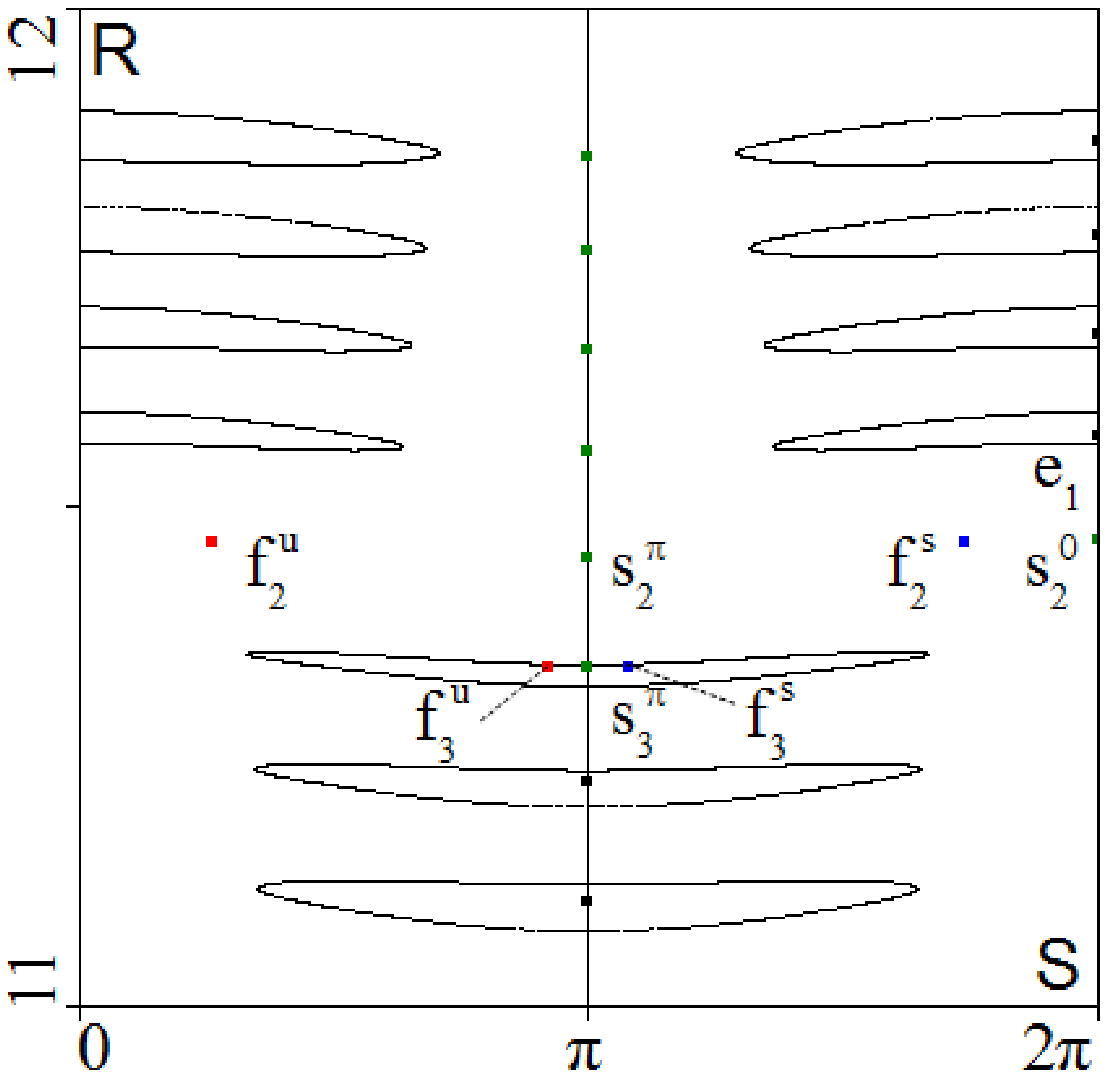} \\ b) {\footnotesize $\varepsilon = 0.015$}}
\end{minipage}
\vfill
\begin{minipage}[h]{0.49\linewidth}
\center{\includegraphics[width=1\linewidth]{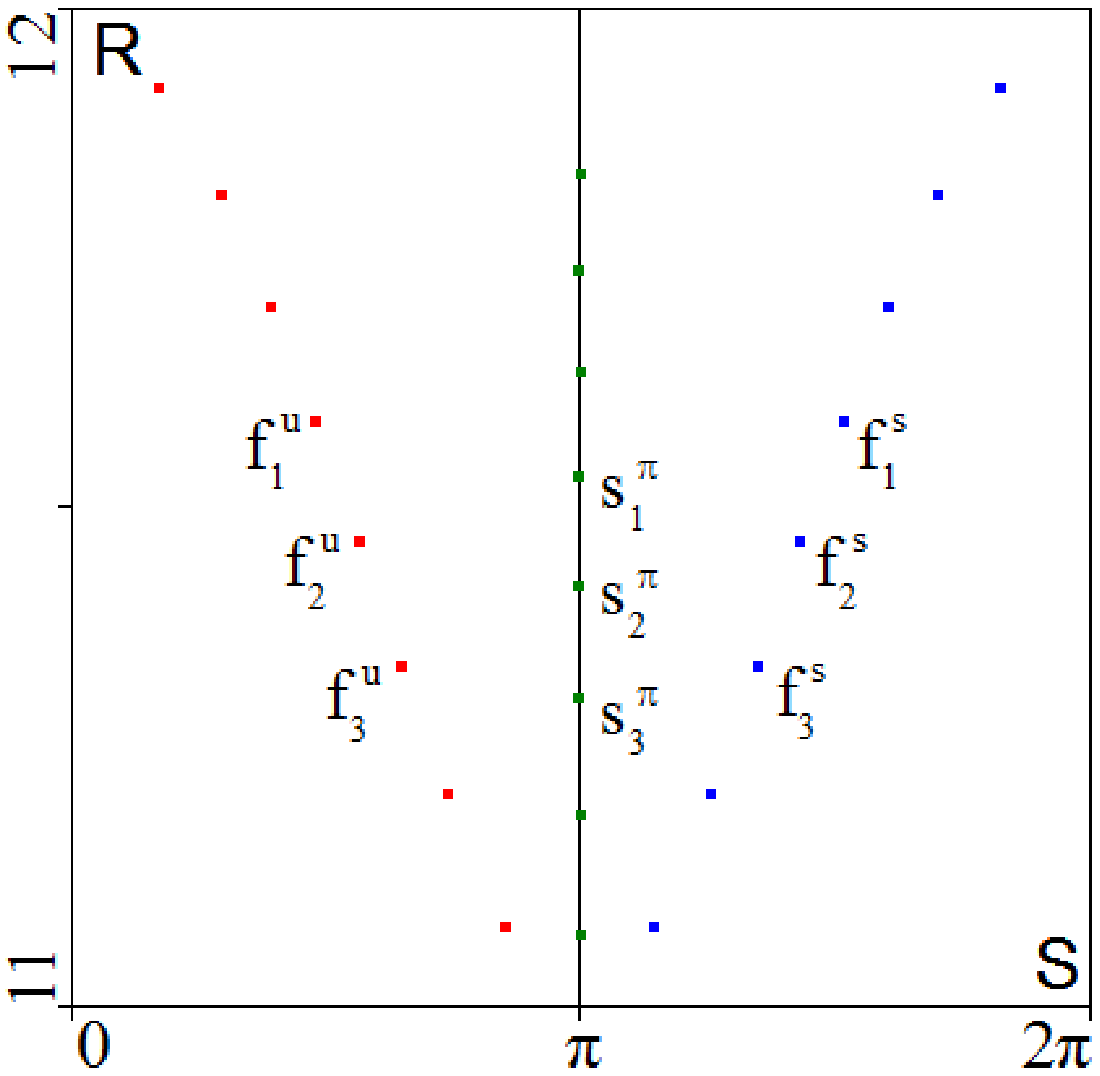} \\ c) {\footnotesize $\varepsilon = 0.1$}}
\end{minipage}
\hfill
\begin{minipage}[h]{0.49\linewidth}
\center{\includegraphics[width=1\linewidth]{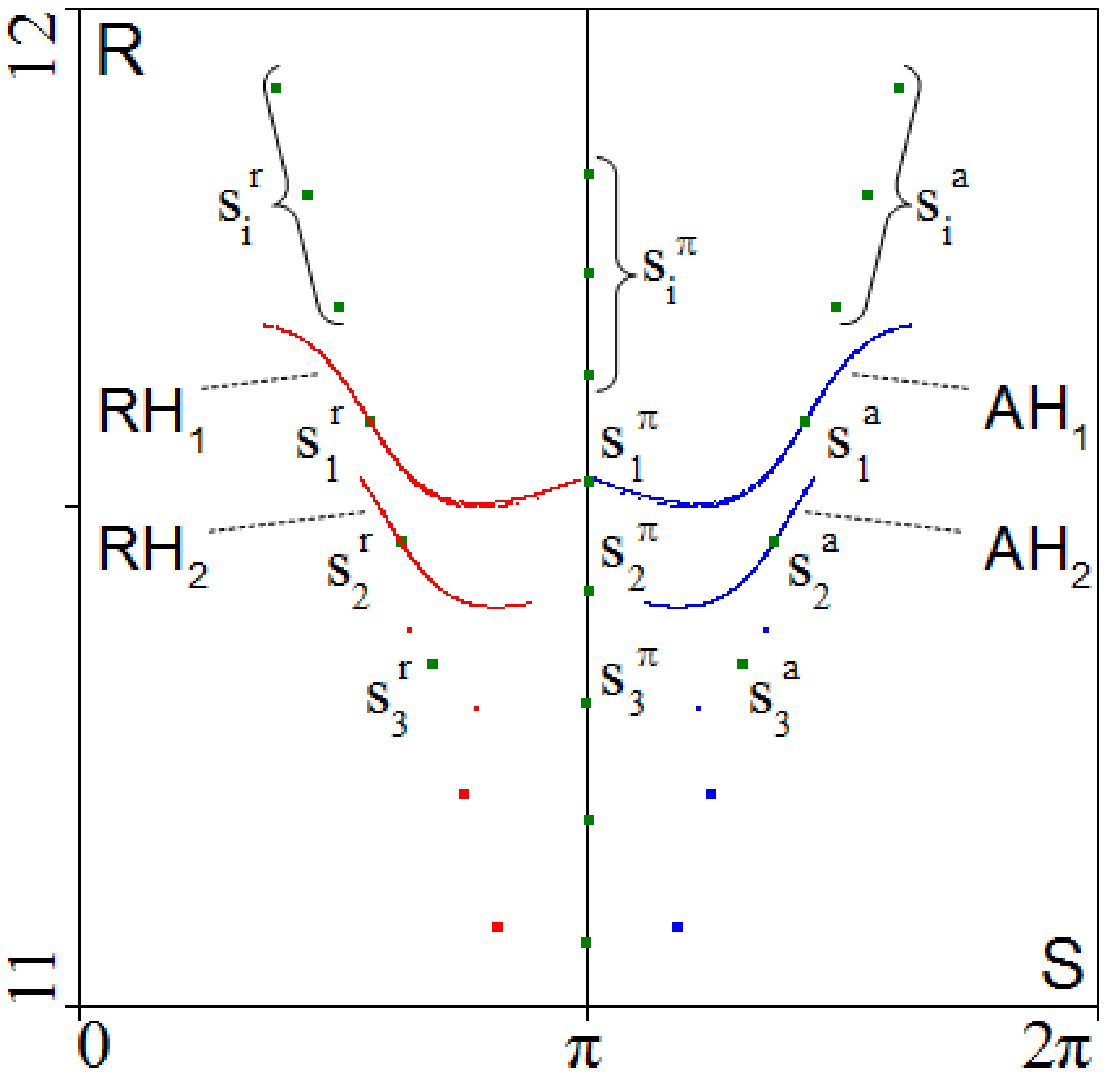} \\ c) {\footnotesize $\varepsilon = 0.146$}}
\end{minipage}
\caption{{\footnotesize The phase spaces of the system \eqref{eq:mainEq} for different $\varepsilon$.
(a) Resonant tori for small $\varepsilon$; (b) elliptic points $e_2$ and $e_3$ undergo symmetry breaking bifurcations due to which symmetric saddles $s_2^{0}$, $s_3^{\pi}$ and also pairs ($f_2^s, f_2^u$) and ($f_3^s, f_3^u$) of stable and completely unstable points appear; (c) 8 stable $f_i^s$ and 8 completely unstable $f_i^u$ fixed points coexist in the system; (d) Henon-like attractors $AH_1$ and $AH_2$ (together with Henon-like repellers $RH_1$ and $RH_2$) born from stable points $f_1^s$  $f_2^s$ (completely unstable points $f_1^u$  $f_2^u$).}}
\label{Fig:IntegrableAndNot}
\end{figure}

For $\varepsilon = 0$ the system \eqref{eq:mainEq} describing the motion of unperturbed vortices is integrable \cite{Yates77} and its phase space foliated into invariant tori. When $\varepsilon > 0$ some tori become resonance and the pair of symmetric saddle and elliptic orbits appear, see Fig. \ref{Fig:IntegrableAndNot}a. With a further increase in the parameter $\varepsilon$ elliptic points $e_i$ undergo symmetry breaking bifurcations due to which these points become symmetric saddles and in their neighborhood stable $f^s_i$ and completely unstable $f^u_i$ fixed points appear (see, Fig. \ref{Fig:IntegrableAndNot}b). Note, that for sufficiently large values of the parameter $\varepsilon$ there are coexist 8 stable $f^s_i$ (and also 8 completely unstable $f^u_i$) fixed points, see Fig. \ref{Fig:IntegrableAndNot}c.

With increasing of $\varepsilon$ Henon-like attractors $AH_i$ and Henon-like repellers $RH_i$ are born from points $f^s_i$ and $f^u_i$ due to scenario presented in Sec. \ref{sec:Scenario}. Figure \ref{Fig:IntegrableAndNot}d shows the coexisting of two Henon-like attractors $AH_1$ and $AH_2$ with stable fixed (and periodic) points for $\varepsilon = 0.146$. In this Figure $s_i^a$ and $s_i^r$ -- saddle point which appear after period doubling bifurcation of stable $f_i^s$ and completely unstable points $f_i^u$, respectively. We note, that saddle points $s_i^a$ have Jacobian less than 1, while saddle points $s_i^r$ have Jacobian greater than one.

When $\varepsilon = \varepsilon_{cris1} \approx 0.14635$ the attractor $AH_1$ and the repeller $RH_1$ undergo crisis due to which these two sets begin to intersect. 
It is important to note, that the intersections of attractors and repellers in the system appear due to homoclinic bifurcations. Further we will describe such bifurcations in details, but first, let us recall two well-known facts related to an evolution of homoclinic attractors:
\begin{itemize}
\item in many cases the boundary of basins of attraction for homoclinic attractors is formed by stable manifolds of some saddle points;
\item a homoclinic attractor is the closure of the unstable manifold of one of its saddle points.
\end{itemize}

\begin{remark}
Both these properties hold for homoclinic Henon-like attractors. In particular, in the Henon map $\bar x = y,\bar y = M - bx - y^2$ with $|b| < 1$ the exact boundary of absorbing domain for the attractor is composed by one (or two, if $b < 0$) of stable manifolds of the saddle fixed point $O_1$ having a positive unstable multiplier, while the homoclinic attractor contains other saddle fixed point $O_2$ with negative unstable multiplier. Here the homoclinic Henon-like attractor $AH$ can be defined as a
the prolongation of the point $O_2$, what means that $AH$ contains (or coincide with) the closure of the unstable manifold $W^U(O_2)$\footnote{Recall that the prolongation of a set $Q$ is a closed, invariant stable set consisting of all points attainable from $Q$ by $\varepsilon^+$-orbits for arbitrarily small $\varepsilon$ \cite{GonTur2017}. Here $\varepsilon^+$-orbits are $\varepsilon$-orbits for forward iteration of the map; $\varepsilon > 0$.}.
\end{remark}

\begin{figure}[!ht]
\begin{minipage}[h]{0.49\linewidth}
\center{\includegraphics[width=1\linewidth]{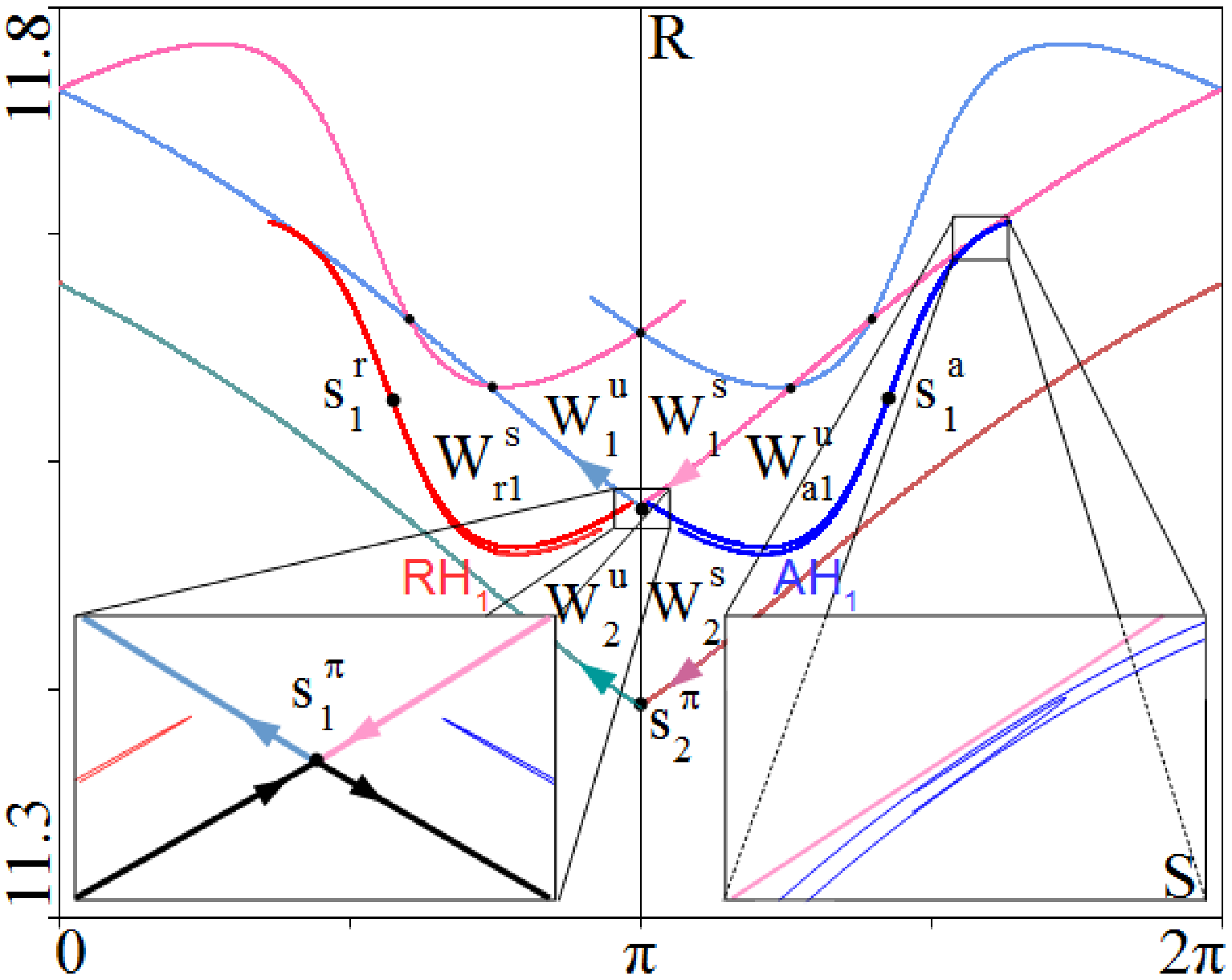} \\ b) {\footnotesize $\varepsilon = 0.1463$}}
\end{minipage}
\hfill
\begin{minipage}[h]{0.49\linewidth}
\center{\includegraphics[width=1\linewidth]{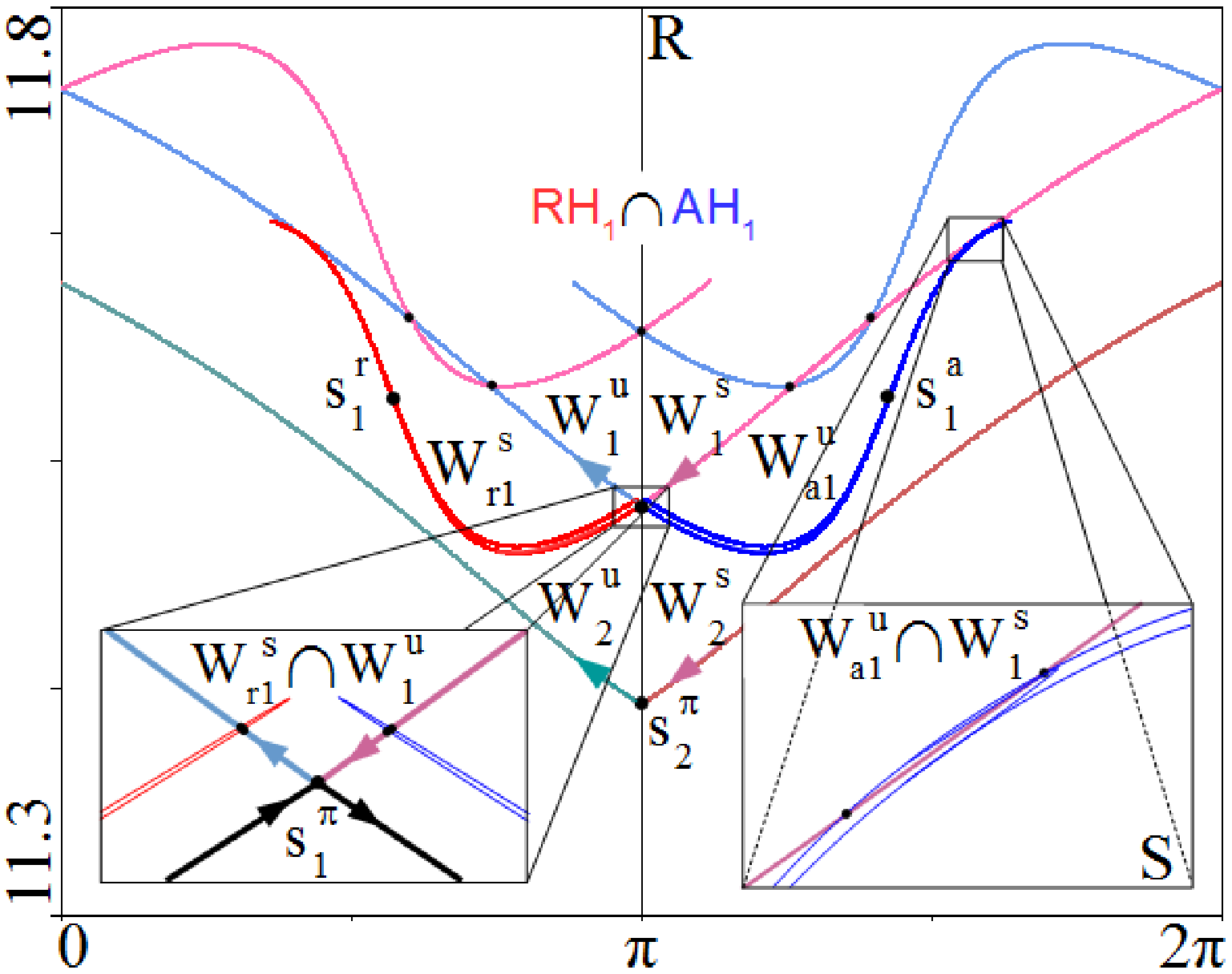} \\ d) {\footnotesize $\varepsilon = 0.1464$}}
\end{minipage}
\caption{{\footnotesize The location of manifolds of the saddle points $s_1^{\pi}, s_2^{\pi}, s_1^a$ and $s_1^r$. $W^u_{a1}$ is the unstable manifold of the saddle point $s_1^a$ belonging to the attractor $AH_1$; $W^s_1$ and $W^s_2$ are the stable manifolds forming the boundary of the basin of attraction for this attractor. $W^s_{r1}$ is the stable manifold of $s_1^r$ belonging to the repeller $RH_1$, while $W^u_1$ and $W^u_2$ are unstable manifolds forming the ``basin of repulsion'' for $RH_1$.
}}
\label{fig:MD_scenario1}
\end{figure}

In the case under consideration the attractor $AH_1$ belongs to the closure of the unstable manifold $W^u_{a1}$ of the saddle fixed point $s_1^a$. The basin of attraction of $AH_1$ is bounded from above by the stable manifold $W^s_1$ of the symmetric saddle point $s_1^{\pi}$, and from below -- by the stable manifold $W^s_2$ of the saddle point $s_2^{\pi}$ (see Fig. \ref{fig:MD_scenario1}a). We also note that the manifolds $W^s_1$ and $W^u_1$ intersect transversally.

\begin{figure}[!ht]
\begin{minipage}[h]{0.49\linewidth}
\center{\includegraphics[width=1\linewidth]{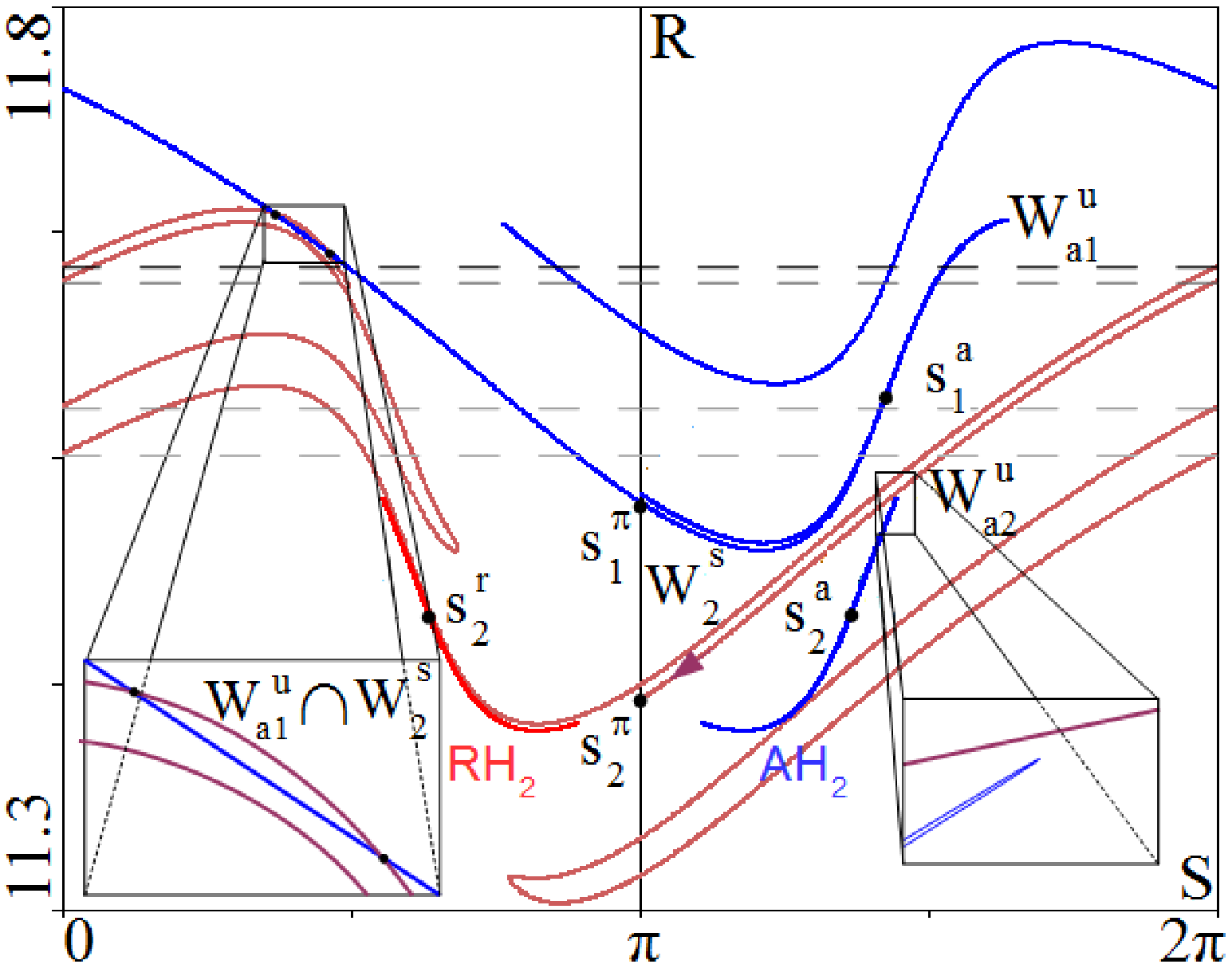} \\ b) {\footnotesize $\varepsilon = 0.1465$}}
\end{minipage}
\hfill
\begin{minipage}[h]{0.49\linewidth}
\center{\includegraphics[width=1\linewidth]{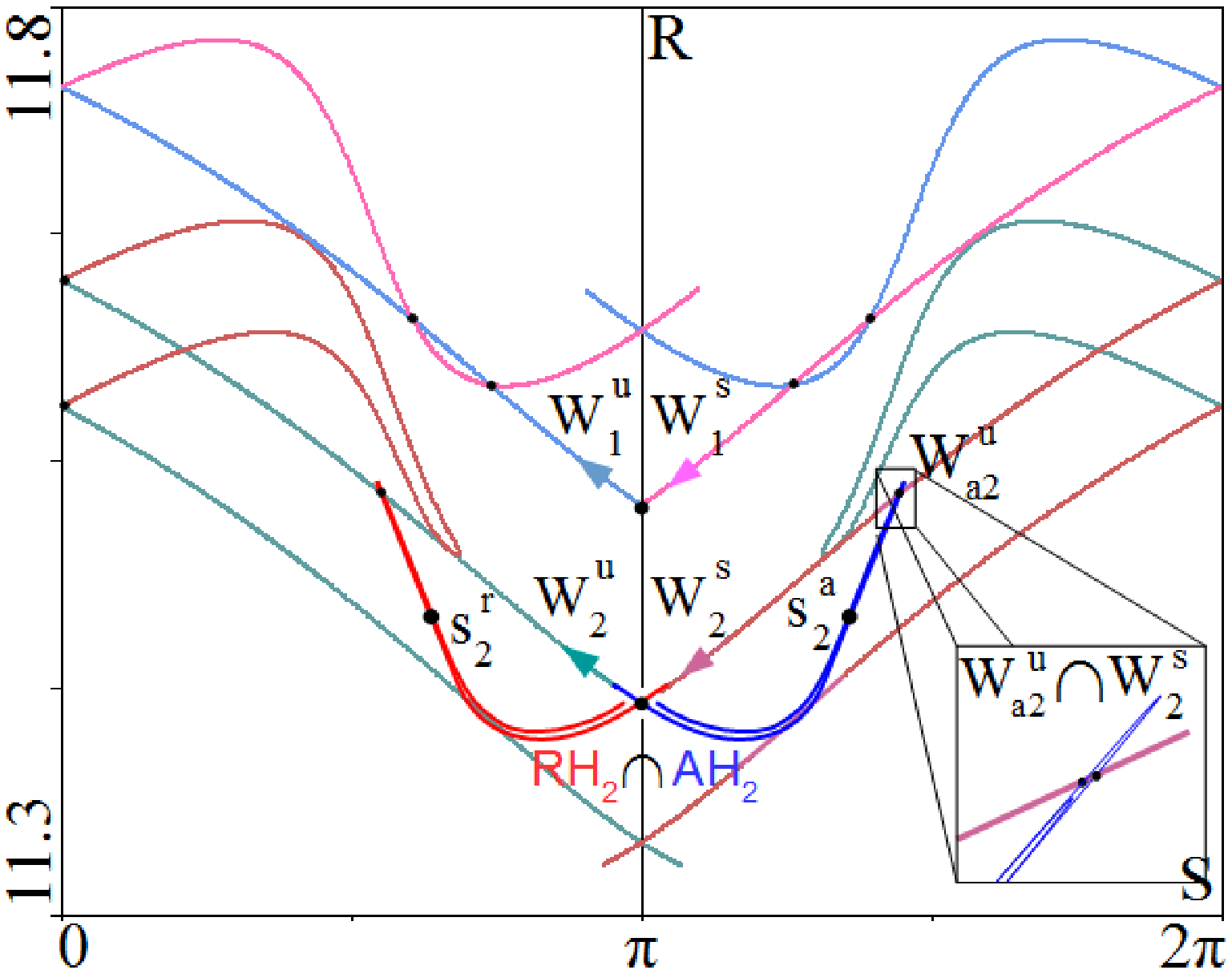} \\ d) {\footnotesize $\varepsilon = 0.1482$}}
\end{minipage}
\caption{{\footnotesize The location of manifold of the saddle points $s_1^{\pi}, s_2^{\pi}, s_1^a$ and $s_1^r$. $W^u_{a1}$ is the unstable manifold of the saddle point $s_1^a$ belonging to the attractor $AH_1$; $W^s_1$ and $W^s_2$ are the stable manifolds forming the boundary of the basin of attraction for this attractor. $W^s_{r1}$ is the stable manifold of $s_1^r$ belonging to the repeller $RH_1$, while $W^u_1$ and $W^u_2$ are unstable manifolds forming the ``basin of repulsion'' for $RH_1$.
}}
\label{fig:MD_scenario2}
\end{figure}

When $\varepsilon > \varepsilon_{cris1}$, intersections between $W^s_1$ and $W^u_{a1}$, as well as $W^u_1$ and $W^s_{r1}$ appear (see Fig. \ref{fig:MD_scenario1}b when $\varepsilon = 0.1464$) and, as a result, the attractor $AH_1$ collides with the upper boundary of its basin of attraction while the repeller $RH_1$ collides with the upper boundary of its ``basin of repulsion''. We note that after this collision the transversal intersection between $W^s_1$ and $W^u_1$ remains and thus we can state here that the attractor $AH_1$ intersects with the repeller $RH_1$.
But in this case such an intersection is not visible due to the existing of homoclinic Henon-like attractor $AH_2$ which attracts almost all trajectories from the neighborhood of $AH_1$ in forward time and $RH_2$ which attracts trajectories in backward time (here after the collision of $AH_1$ and $RH_1$ the unstable manifold $W^u_{a1}$ also intersects with the stable manifold $W^s_2$ of the symmetric saddle point $s_2^{\pi}$, see Fig. \ref{fig:MD_scenario2}a, and this intersection gives a transition mechanism from $AH_1$ to $AH_2$).

When $\varepsilon = \varepsilon_{cris2}$ the unstable manifold $W_{a2}^u$ of $s_2^a$, belonging to the attractor $AH_2$ touches the stable manifold $W_2^s$ which forms the boundary of the basin of attraction for $AH_2$. The same non-transversal heteroclinic tangency appears between another pair of manifolds: $W_{r2}^s$ and $W_2^u$, one of which compose the repeller $RH_2$, while another forms the boundary of the basin for $RH_2$. Thus, for $\varepsilon > \varepsilon_{cris2}$ the attractor $AH_2$ intersects with the repeller $RH_2$, and moreover, due to homoclinic tangle, these two sets intersect with $AH_1 \cap RH_2$ and we have an intersection of two attractors ($AH_1$ and $AH_2$) with two repellers ($RH_1$ and $RH_2$).

But this intersection also is not visible due to existing of other Henon-like attractor $AH_3$ and Henon-like repeller $RH_3$, which appear from $f_3^s$ and $f_3^u$ and, in its turn, collide by the same way as $AH_1$, $RH_1$, $AH_2$ and $RH_2$, giving more complex global connection between different homoclinic attractors and repellers, see Fig. \ref{fig:MD_Vortex}b.

\begin{figure}[!ht]
\begin{minipage}[h]{0.49\linewidth}
\center{\includegraphics[width=1\linewidth]{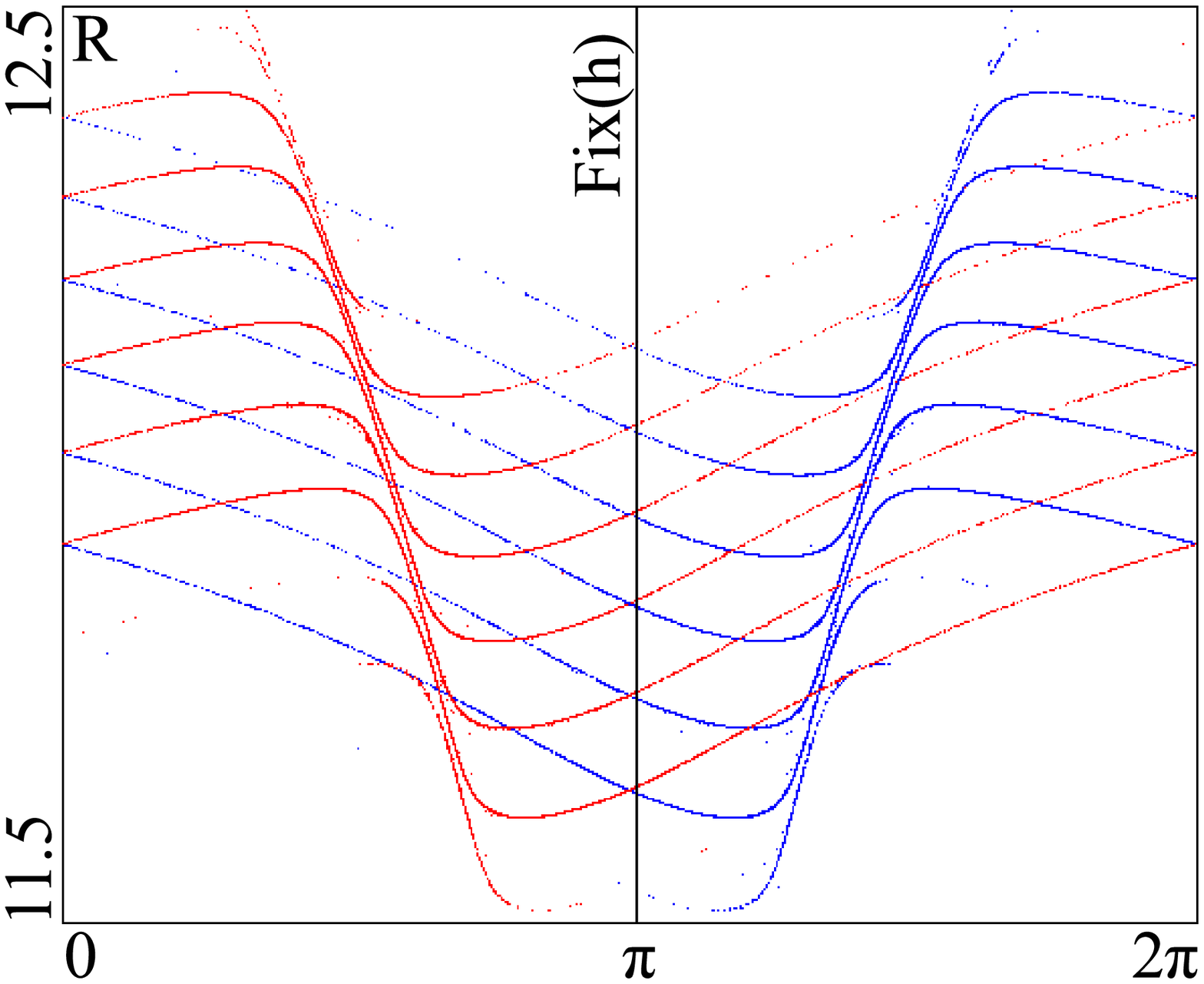} \\ b) {\footnotesize $\varepsilon = 0.23$}}
\end{minipage}
\hfill
\begin{minipage}[h]{0.49\linewidth}
\center{\includegraphics[width=1\linewidth]{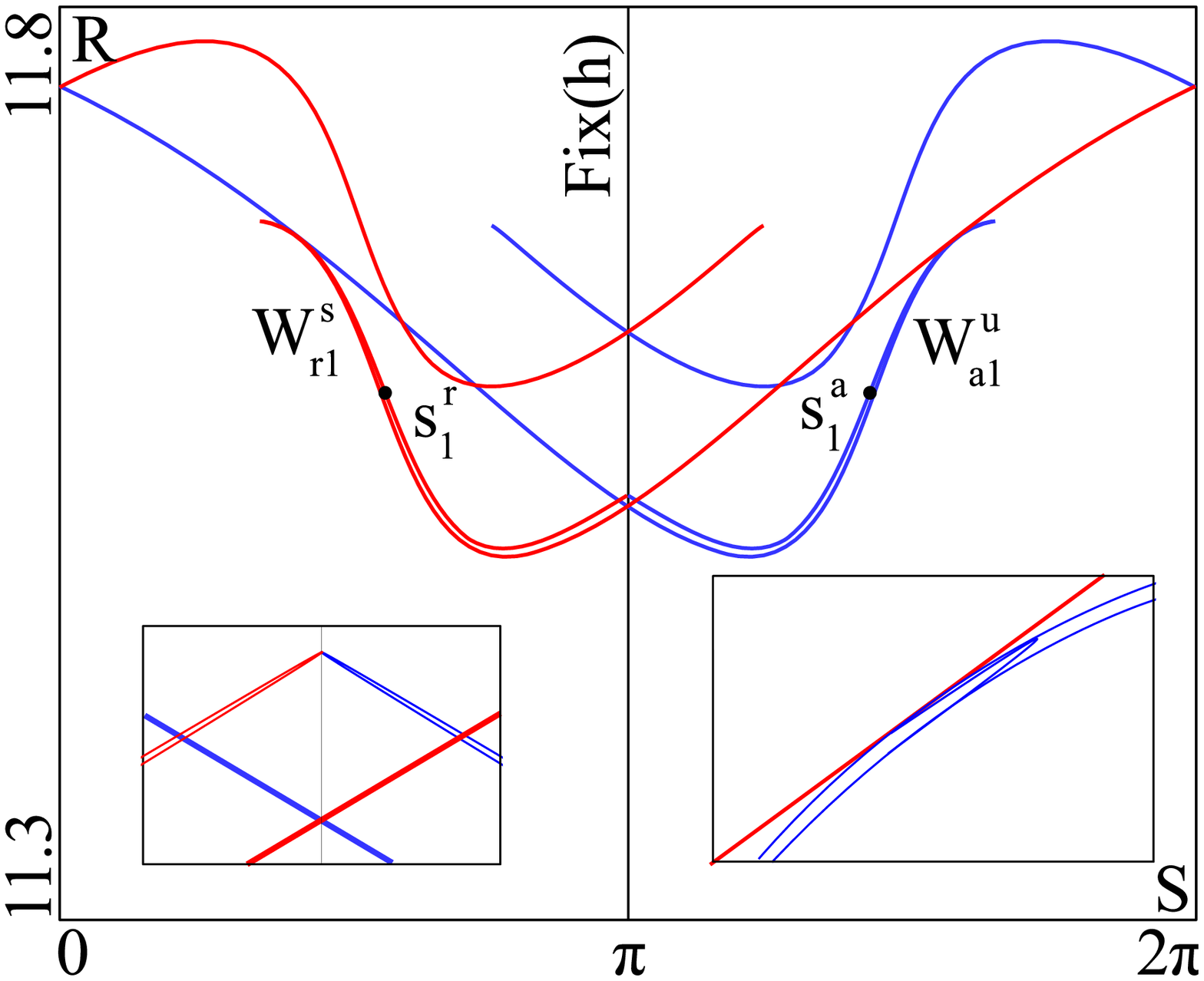} \\ d) {\footnotesize $\varepsilon = 0.14655$}}
\end{minipage}
\caption{{(a) Mixed dynamics after the collion of all homoclinic Henon-like attractors $AH_i$ and all Henon-like repellers $RH_i$; (b) non-transversal heteroclinic cycle containing a pair of non-conservative saddles $s_1^a$, $s_1^r$ with $J < 1$ and $J > 1$, respectively.
}}
\label{fig:MD_final}
\end{figure}

Such a joining of attractors and repellers ends at $\varepsilon > \varepsilon_{cris_8} \approx 0.206$, when all 8 attractors $AH_i$ and repellers $RH_i$ undergo crisis and form a complex set consisting of the intersection of all these attractors and repellers. Since there are no other attractors in the neighborhood of this intersection, the mixed dynamics becomes visible\footnote{It is possible, that mixed dynamics manifests itself before the collision of all 8 attractors and repellers. For example, Fig. \ref{fig:MD_Vortex}b shows mixed dynamics which we obtain with 50000 forward and backward iterations after collision of 3 attractors $AH_i, i = 1, 2, 3$ and 3 repellers $AH_i, i = 1, 2, 3$. In this case it is difficult to understand either we have a coexisting of mixed dynamics with other (in particular, Henon-like) attractors and repeller or mixed dynamics here is transient chaos and finally trajectories evolve to other attractor.}, see Fig. \ref{fig:MD_final}a.

Finally we note, that after each attractor-repeller collision heteroclinic cycles, containing a pair of non-conservative saddles $s_i^a$, $s_i^r$ with $J < 1$ and $J > 1$, appears (see, for example, Fig. \ref{fig:MD_final}b, after the collision of $AH_1$ and $RH_1$, when $\varepsilon = 0.14655$). In the accordance with theorems from \cite{LambStenkin2004, GonDelsh2012} the bifurcations of such cycles lead to a simultaneous birth of infinitely many periodic attractors, repellers and elliptic orbits. Moreover, the closures of these sets have non-empty intersection. From the other hand, at least one area-preserving saddle point $s_i^{\pi}$ (or $s_i^{0}$) belongs to the intersection of $AH_i$ and $RH_i$ after their collision. As it is known from the paper \cite{GonTur2017}, such points belong to a {\it reversible core}, which, by theorem 1 from this paper, contains a limit of an infinite sequence of attractors and repellers. So this theorem gives us another evidence of infiniteness and inseparability of attractors and repellers in the system.


\section{Conclusion}

In this paper, a phenomenological scenario of the emergence of mixed dynamics of a new, strongly dissipative, type is proposed for two-dimensional reversible maps. For this type of mixed dynamics the numerically obtatined chaotic attractor and the chaotic repeller are very different from each other, and their intersection forms a ``thin'' set. The implementation of such a scenario is demonstrated on a system describing the motion of two point vortices perturbed by a wave and external shear flow. It is shown that mixed dynamics in this system contain both conservative and dissipative periodic orbits, while closures of orbits of different types have non-empty intersection. Moreover, in contrast to previously observed types of mixed dynamics, the Jacobians of dissipative orbits in this case can be essentially differ from 1. Thus, it can be argued that strongly dissipative mixed dynamics obviously can not be regarded as quasi-conservative chaos, but represents a new, third (along with conservative and dissipative) type of chaos. Moreover, we assume that this type of chaotic behavior is typical for strongly dissipative reversible systems.

{\bf Acknowledgements}

The author thanks S.V. Gonchenko for valuable advices and comments and E.V. Vetchanin for providing the system in which the phenomenon of the emergence of mixed dynamics was demonstrated.

This paper was supported by grant 17-11-01041 of the RSF. Numerical experiments in Section 4 were supported by RFBR grant 16-01-00364. Also the author was partially supported by the Basic Research Program at the
National Research University Higher School of Economics, project 90 in 2017, and by the Dynasty Foundation.


\begin{thebibliography}{99}{}

\bibitem{Conley78}
Conley C. C.
Isolated invariant sets and the Morse index // American Mathematical Soc., 1978, No. 38.
%
\bibitem{Ruelle1981}
Ruelle D.
Small random perturbations of dynamical systems and the definition of attractors //
Communications in Mathematical Physics. – 1981. – Vol. 82. – No. 1. – pp. 137-151.
%
\bibitem{Hurley1982}
Hurley M.
Attractors: persistence, and density of their basins //
Transactions of the American Mathematical Society. – 1982. – Vol. 269. – No. 1. – pp. 247-271.
%
\bibitem{GonTur2017}
Gonchenko S., Turaev D.
On three types of dynamics, and the notion of attractor //
Trudy Matematicheskogo Instituta imeni V.A. Steklova, 2017, Vol. 297, pp. 133--157.
%
\bibitem{GonShilSten2002}
Gonchenko S. V., Shilnikov L. P., Stenkin O. V.
On Newhouse regions with infinitely many stable and unstable invariant tori //
Proceedings of the Int. Conf. ``Progress in Nonlinear Science'', 2002, Vol., No., pp. 2-6.
%
\bibitem{GonShilSten2006}
Gonchenko S. V., Stenkin O. V., Shilnikov L. P.
On the existence of infinitely many stable and unstable invariant tori for systems from Newhouse regions with heteroclinic tangencies //
Nelineinaya Dinamika [Russian Journal of Nonlinear Dynamics]. – 2006. – Vol. 2. – No. 1. – pp. 3-25.
%
\bibitem{GonDelsh2012}
Delshams A. et al.
Abundance of attracting, repelling and elliptic periodic orbits in two-dimensional reversible maps //
Nonlinearity. – 2012. – Vol. 26. – No. 1. – pp. 1.
%
\bibitem{GonShilTur97}
Gonchenko S. V., Shilnikov L. P., Turaev D. V.
On Newhouse domains of two-dimensional diffeomorphisms which are close to a diffeomorphism with a structurally unstable heteroclinic cycle //
Proc. Steklov Inst. Math. – 1997. – Vol. 216. – pp. 70--118.
%
\bibitem{LambStenkin2004}
Lamb J. S. W., Stenkin O. V.
Newhouse regions for reversible systems with infinitely many stable, unstable and elliptic periodic orbits //
Nonlinearity. – 2004. – Vol. 17. – No. 4. – pp. 1217.
%
\bibitem{PolitiOppoBadii86}
Politi A., Oppo G. L., Badii R.
Coexistence of conservative and dissipative behavior in reversible dynamical systems //
Physical Review A. – 1986. – Vol. 33. – No. 6. – p. 4055.
%
\bibitem{QuispelRoberts92}
Roberts J. A. G., Quispel G. R. W.
Chaos and time-reversal symmetry. Order and chaos in reversible dynamical systems //
Physics Reports. – 1992. – Vol. 216. – No. 2-3. – pp. 63--177.
%
\bibitem{LambRoberts98}
Lamb J. S. W., Roberts J. A. G.
Time-reversal symmetry in dynamical systems: a survey //
Physica D: Nonlinear Phenomena. – 1998. – Vol. 112. – No. 1. – pp. 1--39.
%
\bibitem{PikTop2002}
Topaj D., Pikovsky A.
Reversibility vs. synchronization in oscillator lattices //
Physica D: Nonlinear Phenomena. – 2002. – Vol. 170. – No. 2. – pp. 118--130.
%
\bibitem{GonGonKazTur2017}
Gonchenko A.S., Gonchenko S.V., Kazakov A.O., Turaev D.V.
On the phenomenon of mixed dynamics in Pikovsky–Topaj system of coupled rotators //
Physica D: Nonlinear Phenomena. – 2017. – Vol. 350. – pp. 45--57.
%
\bibitem{GonGonKaz2013}
Gonchenko A. S., Gonchenko S. V., Kazakov A. O.
Richness of chaotic dynamics in nonholonomic models of a Celtic stone //
Regular and Chaotic Dynamics. – 2013. – Vol. 18. – No. 5. – pp. 521--538.
%
\bibitem{Kaz2013}
Kazakov A. O.
Strange attractors and mixed dynamics in the problem of an unbalanced rubber ball rolling on a plane //
Regular and Chaotic Dynamics. – 2013. – Vol. 18. – No. 5. – pp. 508--520.
%
\bibitem{Kazakov2015}
Kazakov A.
On chaotic dynamics in the Suslov problem //
Dynamics, bifurcations and chaos 2015 (DBC II): Ext. Abstr. Int. Conf.–Sch., Nizhni Novgorod, July 20–24, 2015. Nizhni Novgorod: Lobachevsky State Univ., 2015. P. 21–30.
%
\bibitem{BizBorKaz2015}
Bizyaev I. A., Borisov A. V., Kazakov A. O.
Dynamics of the suslov problem in a gravitational field: Reversal and strange attractors //
Regular and Chaotic Dynamics. – 2015. – Vol. 20. – No. 5. – pp. 605--626.
%
\bibitem{Vet2017}
Vetchanin  E.V.,  Mamaev I.S. 
Dynamics of two point vortices in an external compressible shear flow // 
Regular and Chaotic Dynamics. - 2017. - Vol. 22. - No. 8, - pp. 893–-908.
%
\bibitem{LermanTuraev2012}
Lerman L. M., Turaev D.
Breakdown of symmetry in reversible systems //
Regular and Chaotic Dynamics. – 2012. – Vol. 17. – No. 3--4. – pp. 318-336.
%
\bibitem{Feigenbaum83}
Feigenbaum M. J.
Universal behavior in nonlinear systems //
Physica D: Nonlinear Phenomena. – 1983. – Vol. 7. – No. 1-3. – pp. 16--39.
%
\bibitem{Henon76}
Henon M.
A two-dimensional mapping with a strange attractor //
The Theory of Chaotic Attractors. – Springer New York, 1976. – pp. 94--102.
%
\bibitem{GonGon2016}
Gonchenko A. S., Gonchenko S. V.
Variety of strange pseudohyperbolic attractors in three-dimensional generalized Henon maps //
Physica D: Nonlinear Phenomena. – 2016. – Vol. 337. – pp. 43-57.
%
\bibitem{GoncharOstapchukTur91}
Gonchar V.Y. , Ostapchuk P.N., Tur A.V., Yanovsky V.V.
Dynamics and stochasticity in a reversible system describing interaction of point vortices with a potential wave //Physics Letters A. – 1991. – Vol. 152. – No. 5-6. – pp. 287-292.
%
\bibitem{KazVet2016}
Vetchanin E. V., Kazakov A. O.
Bifurcations and chaos in the dynamics of two point vortices in an acoustic wave //
International Journal of Bifurcation and Chaos. – 2016. – Vol. 26. – No. 04. – p. 1650063.
%
\bibitem{Yates77}
Yates J. E.
Interaction with and production of sound by vortex flows //
American Institute of Aeronautics and Astronautics Conference. – 1977.



\end{thebibliography}
\end{document}